\documentclass[11pt]{article}
\usepackage[dvips]{graphics}
\usepackage{multirow}
\usepackage{array}
\usepackage{epsfig,rotating}
\usepackage{amssymb,amsmath}
\usepackage{latexsym}
\usepackage[usenames]{color}

 \newtheorem{theorem}{\textbf{Theorem}}[section]
 
 \newtheorem{lemma}{\textbf{Lemma}}[section]

 \newtheorem{exm}{\textbf{Example}}

 \newtheorem{rem}{\textbf{Remark}}

\numberwithin{equation}{section}

\newcommand{\be}{\begin{equation}}
\newcommand{\ee}{\end{equation}}
\newcommand{\beaa}{\begin{eqnarray*}}
\newcommand{\eeaa}{\end{eqnarray*}}
\newcommand{\bea}{\begin{eqnarray}}
\newcommand{\eea}{\end{eqnarray}}
\newcommand{\lbl}{\label}

\newcommand{\bei}{\begin{itemize}}
\newcommand{\eei}{\end{itemize}}

\def\Var{\mathrm{Var}}

\newcommand{\bz}{\mathbf{z}}
\begin{document}
\thispagestyle{empty}

\noindent{\bf \Large Limiting Empirical Spectral Distribution for
Products \\of Rectangular Matrices}

\vspace{10pt} \noindent \textbf{Yongcheng Qi$^1$, ~~~Hongru
Zhao$^2$}

\vspace{20pt}

{\footnotesize
\noindent $^1$Department of Mathematics and
Statistics, University of Minnesota Duluth, 1117 University Drive,
Duluth, MN 55812, USA. Email: yqi@d.umn.edu  (corresponding author)

\noindent $^2$Department of Mathematics and Statistics, University
of Minnesota Duluth, 1117 University Drive, MN 55812, USA. Email:
zhao1118@d.umn.edu.
}

\date{\today}

\vspace{20pt}




{\small

\noindent{\bf Abstract.} In this paper, we consider $m$ independent
random rectangular matrices whose entries are independent and
identically distributed standard complex Gaussian random variables
and assume the product of the $m$ rectangular matrices is an $n$ by
$n$ square matrix. We study the limiting empirical spectral
distributions of the product where the dimension of the product
matrix goes to infinity, and $m$ may change with the dimension of
the product matrix and diverge. We give a complete description for
the limiting distribution of the empirical spectral distributions
for the product matrix and illustrate some examples.

\vspace{10pt}

 \noindent \textbf{Keywords:\/} Empirical spectral distribution,
Eigenvalues, Product of rectangular matrices, Non-Hermitian random
matrix

\vspace{10pt}

\noindent\textbf{AMS 2010 Subject Classification: \/} 15B52, 60F99, 60G57 \\
}

\newpage

\section{Introduction}\label{intro}

The study of Random Matrix Theory was initialized by
Wishart~\cite{Wishart} for statistical analysis of large samples.
Wigner~\cite{Wigner} found applications for random Hermitian matrix
in nuclear physics. Subsequential applications include condensed
matter physics (Beenakker~\cite{Been1997}), number theory (Mezzadri
and Snaith~\cite{MS2005}), wireless communications (Couillet and
Debbah~\cite{CD}, and high dimensional statistics
(Johnstone~\cite{John2001, John2008}, Jiang~\cite{Jiang09}), quantum
chromodynamics, chaotic quantum systems and growth processes (see,
e.g., Akemann, Baik and Francesco~\cite{ABF2011}).

There are two major directions for the study for random matrices,
including the empirical spectral distributions and the spectral
radii. The classical semi-circular law was first introduced by
Wigner, and then Ginibre~\cite{Ginibre} established the circle law
for Ginibre ensembles. Since then, the assumptions were relaxed
subsequently in the papers by Girko~\cite{Girko},  Bai~\cite{Bai},
Pan and Zhou~\cite{Pan}, and G\"otze and Tikhomirov~\cite{GF-TA}.
Tao and Vu~\cite{Tao} proved the circular law under the second
moment condition. For the spectral radii, Tracy and Widom
established the so-called Tracy-Widom laws for the limiting
distributions for the three Hermitian matrices (Gaussian orthogonal
ensemble, Gaussian unitary ensemble and Gaussian symplectic
ensemble); see Tracy and Widom~\cite{Tracy94, Tracy96}. Other work
in this direct includes Rider~\cite{Rider2003, Rider2004} and Rider
and Sinclair~\cite{RS2014}.

Products of random matrices are particularly of interest in recent
research. Ipsen~\cite{Ipsen} provided several applications,
including wireless telecommunication, disordered spin chain, the
stability of large complex system, quantum transport in disordered
wires, symplectic maps and Hamiltonian mechanics, quantum
chromo-dynamics at non-zero chemical potential. G\"{o}tze and
Tikhomirov~\cite{Goetz}, Bordenave~\cite{Bor}, O'Rourke and
Soshnikov~\cite{Rourke} and O'Rourke {\it et al.}~\cite{Rourke14}
found the limiting empirical spectral distribution for the product
from the complex Ginibre ensemble when $m$ is fixed.  Two recent
papers by Jiang and Qi~\cite{JiangQi2017, JiangQi2019} considered
the spectral radii and limiting empirical spectral distribution for
the product of complex Ginibre ensembles and the product of
truncations of independent Haar unitary matrices by allowing $m$ to
change. G\"otze, K\"osters and Tikhomirov~\cite{GKT2015},
Zeng~\cite{Zeng2016}, and Chang and Qi~\cite{ChangQi2017} studied
the limiting empirical distribution of product of the spherical
ensemble. Chang, Li and Qi~\cite{ChangLiQi2018} investigated the
limiting distribution of the spectral radii for product of matrices
from the spherical ensemble.

In this paper, we consider the product of $m$ random rectangular
matrices with independent and identically distributed (i.i.d.)
complex Gaussian entries and investigate the limiting empirical
spectral distributions. Adhikari {\it et al.}~\cite{ARRS2016}
obtained the joint density function for the eigenvalues and found
the limit of the expected empirical distributions when $m$ is a
fixed integer, and Zeng~\cite{Zeng2017} obtained the limiting
empirical spectral distribution. Lambert~\cite{Lambert2018}
established that the empirical distribution for square singular
values converges to certain generalizations of the Fuss-Catalan
distribution and that the maximum of the square singular values
converges to the edge point of the Fuss-Catalan distribution. Very
recently, Qi and Xie~\cite{QiXie2019} obtained the limiting
distributions for spectral radii for products of rectangular
matrices when $m$ changes with the dimension of the product
matrices.


The rest of the paper is organized as follows.  In
Section~\ref{main},  we introduce empirical spectral distributions
for scaled eigenvalues from the production of independent random
rectangular matrices and present a general result on the convergence
of the empirical spectral distributions. We further investigate the
limiting distributions and obtain all types of distributions and
provide conditions when these distributions can be obtained.  We
also give a few illustrative examples.  Proofs for the main results
are given in Section~\ref{proofs}.

\section{Main Results}\label{main}

In this paper, we consider $m$ independent rectangular matrices,
$\mathbf{X}_j,~ 1 \leq j \leq m$, namely $\mathbf{X}_j$ is an
$n_j\times n_{j+1}$ matrix for $1\le j\le m$, where $n_1, \cdots,
n_{m+1}$ are positive integers, and all entries of the $m$ matrices
are independent and identically distributed standard complex normal
random variables. We assume $n_1=n_{m+1}=:n$ so that the product
\[
\mathbf{X}^{(m)}=\prod^m_{j=1}\mathbf{X}_j\]
 is an $n\times n$ square
matrix. We also assume $n=\min_{1\le j\le m+1}n_j$.  In this case,
the product matrix $\mathbf{X}^{(m)}$ is of full rank.

Denote the $n$ eigenvalues of $\mathbf{X}$ as $\mathbf{z}_{1},
\cdots, \mathbf{z}_{n}$, and set $l_{j}=n_j-n\geq 0$, $j=1, \cdots,
m$. It follows from Theorem 2 of Adhikari {\it et
al.}~\cite{ARRS2016} that the joint density function for
$\mathbf{z}_{1}, \cdots, \mathbf{z}_{n}$ is given by
\begin{equation}\label{model}
 p(z_{1}, \cdots,
z_{n})=C \prod_{1 \leq j<k \leq n}\left|z_{j}-z_{k}\right|^{2}
\prod_{j=1}^{n} w_{m}^{\left(l_{1}, \cdots,
l_{m}\right)}\left(\left|z_{j}\right|\right)
\end{equation}
with respect to Lebesgue measure on $\mathbb{C}^{n}$, where $C$ is a
normalizing constant such that $p(z_{1}, \cdots, z_{n})$ is a
probability density function, and function $w_{m}^{\left(l_{1},
\cdots, l_{m}\right)}(z)$ can be obtained recursively by
\[
w_{k}^{\left(l_{1}, \cdots, l_{k}\right)}(z)=2 \pi \int_{0}^{\infty}
w_{k-1}^{\left(l_{1}, \cdots,
l_{k-1}\right)}\left(\frac{z}{s}\right)
w_{1}^{\left(l_{k}\right)}(s) \frac{d s}{s}, \quad k \geq 2
\]
with initial $w_{1}^{(l)}(z)=\exp \left(-|z|^{2}\right)|z|^{2 l}$
for any $z$ in the complex plane; see Zeng~\cite{Zeng2016}.

Our objective in the paper is to investigate the limiting empirical
spectral distribution of the product ensemble $\mathbf{X}^{(m)}$
when $n$ tend to infinity. We allow $m$ to change with $n$ and
substitute $m_n$ for $m$ from now on to show its dependence on $n$.

The empirical spectral distribution of $\mathbf{X}^{(m)}$ is the
empirical distribution based on the eigenvalues, $\mathbf{z}_{1},
\cdots, \mathbf{z}_{n}$, of $\mathbf{X}^{(m)}$, i.e.,
\begin{equation}\label{mu*}
\mu_n^*=\frac1n\sum^n_{j=1}\delta_{\mathbf{z}_j/a_n},
\end{equation}
where $a_n>0$ is a sequence of normalizing constants.  When $m_n$
diverges with $n$, the magnitude of $\mathbf{z}_j$'s can go to
infinity exponentially or vanish exponentially. In this case, one
may not be able to find a sequence $a_n$ such that the empirical
measure $\mu_n^*$ converges. Instead, we will define empirical
distribution for scaled eigenvalues as in Jiang and
Qi~\cite{JiangQi2019}.

Note that $\{\mathbf{z}_j;\, 1\leq j \leq n\}$ are complex random
variables. Write
 \begin{equation}\label{argument}
\Theta_j=\arg(\mathbf{z}_j)\in [0, 2\pi)~\mbox{ such that }
~\mathbf{z}_j=|\mathbf{z}_j|\cdot e^{i\Theta_j}
\end{equation}
 for $1\le j\le n$. Further, assume that  $Y_1, \cdots, Y_n$ are independent
random variables and $Y_j$ has a density function proportional to
$y^{j-1}w_{m}^{(l_{1}, \cdots, l_{m})}(y)I(y>0)$. Given a sequence
of positive measurable functions $h_n(r), \, n\geq 1$, which are
defined on $(0,\infty)$, we define the empirical measures for scaled
eigenvalues as follows
\begin{equation}\label{mun}
\mu_n=\frac{1}{n}\sum^n_{j=1}\delta_{(\Theta_j,
h_n(|\mathbf{z}_j|))}~~~\mbox{and}~~~
\nu_n=\frac{1}{n}\sum^n_{j=1}\delta_{h_n(Y_j)}.
\end{equation}

We note that the empirical spectral measure $\mu_n^*$ defined in
\eqref{mu*} is the joint distribution for linearly scaled
eigenvalues, which is the joint empirical distribution based on real
parts and imaginary parts for linearly scaled eigenvalues. The
empirical spectral measure $\mu_n$ defined in \eqref{mun} is the
joint distribution for arguments and scaled moduli of eigenvalues.
The transformation $h_n$ which applies to the moduli of eigenvalues
can be any positive function.  With notation in \eqref{argument}, we
can use $(\Theta_j, h_n(\bz_j))$ to form a new complex number
$h_n(|\bz_j|)e^{i\Theta_j}$.  Therefore, we can define the empirical
measure for scaled eigenvalues $h_n(|\bz_j|)e^{i\Theta_j}$ as
follows
\begin{equation}\label{hatmu}
\hat\mu_n=\frac1n\sum^n_{j=1}\delta_{h_n(|\mathbf{z}_j|)e^{i\Theta_j}}.
\end{equation}
We want to menton that two measures $\hat\mu_n$ and $\mu_n^*$ are
the same when $h_n(r)=r/a_n$.

We will see later that the convergence of $\mu_n$ is equivalent to
that of $\nu_n$. In (\ref{mun}), if $h_n$ is linear, that is,
$h_n(r)=r/a_n$, where $\{a_n, n\geq 1\}$ is a sequence of positive
numbers, we denote the empirical measure of $\mathbf{z}_j$'s by
$\mu_n^*$ as in \eqref{mu*}, and accordingly, we denote the
empirical distribution of $Y_j$'s by
\begin{equation}\label{nun*}
\nu_n^*=\frac{1}{n}\sum^n_{j=1}\delta_{Y_j/a_n}.
\end{equation}

We need the following notations as in the paper by Jiang and
Qi~\cite{JiangQi2019}.

\noindent $\bullet$ Any function $g(z)$ of complex variable
$z=x+iy$, $x,y\in\mathbb{R}$ should be interpreted as a bivariate
function of $(x,y)$:  $g(z)=g(x, y)$.

\noindent $\bullet$ We write $\int_Ag(z)\,dz=\int_Ag(x,y)\,dxdy$ for
any measurable set $A\subset \mathbb{C}.$

\noindent $\bullet$  $\mbox{Unif}(A)$ stands for the uniform
distribution on a set $A$.

\noindent $\bullet$ For a sequence of random probability measures
$\{\tau, \tau_n;\, n\geq 1\}$, we write \bea\lbl{my_baby} \tau_n
\rightsquigarrow \tau\ \ \mbox{if \ $\mathbb{P}$($\tau_n$  converges
weakly to $\tau$ as $n\to\infty$)=1}. \eea When $\tau$ is a
non-random probability measure generated by random variable $X$, we
simply write $\tau_n \rightsquigarrow X$. Review the notation
``$\rightsquigarrow$" in (\ref{my_baby}). The symbol
$\mu_1\otimes\mu_2$ represents the product measure of two measures
$\mu_1$ and $\mu_2$.

For determinantal point processes, Jiang and Qi~\cite{JiangQi2019}
have established a general result on convergence of the empirical
spectral distributions; see Lemma~\ref{nonlinear} in
Section~\ref{proofs}.

It follows form Lemma~\ref{nonlinear} that a common feature for
limiting empirical distributions from determinant point processes is
that the angle and radius of the random vector with the liming
distribution are independent and the convergence of empirical
distributions for the eigenvalues is equivalent to the convergence
of the empirical distribution based the radii of the eigenvalues.

Inspired by Jiang and Qi~\cite{JiangQi2019} and
Zeng~\cite{Zeng2017}, we define a sequence of distribution functions
$F_n(x)$ as follows
\begin{eqnarray}\label{Fn}
    F_{n}(x)=\Big(\prod_{j=1}^{m_n} \frac{nx+l_{j} }{n+l_{j}}\Big)^{1 / \gamma_{n}}=
    \Big(\prod_{j=1}^{m_n} (1-\frac{n}{n_j}(1-x))\Big) ^{1 /\gamma_{n}},~~~~  x
\in[0,1],
\end{eqnarray}
where $\{\gamma_n\}$ is a sequence of positive numbers to be
selected so that $F_n$ has a limit. Note that $F_n(x)$ is continuous
and strictly increasing on $[0,1]$ with $F_{n}(0)=0$ and
$F_{n}(1)=1$. It is easy to see that $F_n$ is a distribution
function on $[0,1]$.  We assume $F_n(x)=0$ when $x<0$ and $F_n(x)=1$
when $x>1$.

We will assume that $F_n(x)$ converges weakly to a distribution
function $F(x)$. This limit is closely related to the limiting
empirical spectral distribution of $\mu_n$ and $\mu_n^*$ defined in
\eqref{mun} and \eqref{mu*}.

A cumulative distribution $F$ is a nondecreasing right-continuous
function, and its generalized inverse defined as
\begin{equation}\label{G-inverse}
F^*(u)=\inf\{x:~ F(x)>u\}, ~~~u\in [0,1)
\end{equation}
Define $F^*(u)=0$ for $u<0$ and $F^*(u)=1$ for $u\ge 1$. One can
show that $F^*(u)$ is also a nondecreasing right-continuous function
and therefore, $F^*$ is also a cumulative distribution function.
When $F$ is continuous and strictly increasing,  $F^*$ is the
regular inverse of $F$.

Recall that $F_n$ converges weakly to a distribution $F$ if and only
if $\lim_{n\to\infty}F_n(x)=F(x)$ for every continuity point $x$ of
$F$. A  probability measure $v$ is induced by $F^*$ if
$\nu((-\infty, u])=F^*(u)$ for all $u$.

The main results of the paper are the following Theorems~\ref{typeI}
and \ref{typeII}.

\begin{theorem}\label{typeI}
Let $\{m_n, n\geq 1\}$ be a sequence of positive integers and
$\gamma_n>0$. Assume that, for any positive integer $k$,
    \begin{equation}\label{vip}
     c_k:=\lim_{n\rightarrow\infty}\frac{1}{\gamma_{n}} \sum_{r=1}^{m_n}
\left(\frac{n}{n_{r}}\right)^k\mbox{ exists}
    \end{equation}
    with $c_1\in (0,\infty)$, and $c_k\in [0, c_1]$ is non-increasing in $k\ge
2$. Define a distribution function $F$ as follows
\begin{equation}\label{F-rep}
F(x)=\exp(-\sum^\infty_{k=1}\frac{c_k}{k}(1-x)^k),~~~~x\in (0,1],
\end{equation}
and its generalized inverse, $F^*$, is given in \eqref{G-inverse}.
  Set
$h_n(x)=\frac{1}{a_n}|x|^{2/\gamma_{n}}$ with
$a_n=\prod_{r=1}^{m_n}n_r^{1/\gamma_{n}}$. Then
$\hat\mu_n\rightsquigarrow \mu$, where $\hat\mu_n$ is defined as in
\eqref{hatmu}, and $\mu$ has a density function $
\frac{f^*(|z|)}{2\pi|z|}I(F(0)\le |z|\le 1)$, where $f^*$ is the
density function of $F^*$ and it can be also determined by
$f^*(x)=1/f(F^*(x))$ with $f(x)=F'(x)$, $x\in (0,1]$.
\end{theorem}

\begin{theorem}\label{typeII}
Let $\{m_n, n\geq 1\}$ be a sequence of positive integers and
$\gamma_n>0$. Assume
  \begin{equation}\label{zero}
    \lim_{n\rightarrow\infty}\frac{1}{\gamma_{n}} \sum_{r=1}^{m_n} \frac{n}{n_{r}}=0.
    \end{equation}
Define $h_n(x)=\frac{1}{a_n}|x|^{2/\gamma_{n}}$ with
$a_n=\prod_{r=1}^{m_n}n_r^{1/\gamma_{n}}$. Then
$\hat\mu_n\rightsquigarrow \mathrm{Unif}(|z|=1)$, where $\hat\mu_n$
is defined as in \eqref{hatmu}.
\end{theorem}

Next, we present some general results on the convergence of the
empirical distribution $F_n$. We will investigate the necessary and
sufficient conditions for the weak convergence of $F_n$,
characterize its limiting distribution $F$ and reveal how the
function $F$ is related to the limit of the empirical measures
$\mu_n$. Theorems~\ref{typeI} and \ref{typeII} are the direct
consequences of the following two theorems.

\begin{theorem}\label{general}
    Let $\{m_n\}$ be an arbitrary sequence of positive
    integers and $\{\gamma_n\}$ be a sequence of positive numbers
    such that $F_n$ converges weakly to a probability distribution
    $F$.  Let $F^*$ denote the generalized inverse of $F$ and $\nu$
    be  a probability measure induced by $F^*$. Define $a_n=\prod_{r=1}^{m_n}n_r^{1/\gamma_{n}}$ and $h_n(x)=\frac{1}{a_n}|x|^{2/\gamma_{n}}$
    in \eqref{mun}.  Then we have $\mu_n \rightsquigarrow
\mathrm{Unif}[0,2 \pi)\otimes \nu$ as $n\to
    \infty$.
\end{theorem}

\begin{theorem}\label{type}
    Let $\{m_n\}$ be a sequence of positive integers, and $\gamma_n$ be any sequence of positive numbers.
    If  $F_n(x)$ converges weakly to a distribution function $F(x)$,
then  $F$ is of one of the following three types

\begin{description}

\item [(Type I).] $F(x)$ is continuous on $[0,1]$, and analytic on
$(0,1)$, with $F(0+)\geq0$,
    $F(1)=1$, and the first derivative $f(x)=F'(x)> 0$  for $x\in
(0,1)$;

\item [(Type II).] $F(0-)=0$, $F(x)=1$ for all $x\in [0,1]$;

\item [(Type III).] $F(1)=1$, $F(x)=0$ for all $x\in [0,1)$.
\end{description}
Furthermore, we have
\begin{itemize}
\item[(a).]
    $F_n(x)$ converges weakly to  a Type I distribution if and only if condition
    \eqref{vip} holds; Under condition \eqref{vip}, the limiting distribution $F$ has
a representation given in \eqref{F-rep}.

\item[(b).] $F_n(x)$ converges weakly to a Type II distribution if and only if
\eqref{zero} holds.

\item[(c).]
    $F_n(x)$ converges weakly to a Type III distribution if and only if
    \begin{equation}\label{infinity}
    \lim_{n\rightarrow\infty}\frac{1}{\gamma_{n}} \sum_{r=1}^{m_n}
\frac{n}{n_{r}}=\infty.
    \end{equation}
\end{itemize}
\end{theorem}



\vspace{20pt}

\begin{rem}
From Theorem~\ref{type}, we can draw the following conclusions.

 \noindent \textbf{a.} If $F$ is of type I,
 $F$ is strictly increasing in $[0,1]$ and its generalized inverse
$F^*$ is given by
    \begin{equation}\label{case-a}
    F^{*}(x)=
\left\{
\begin{array}{ll}
{0,} &\mbox{ if } x < F(0), \\
{F^{-1}(x),} & \mbox{ if } x \in [F(0),1), \\
{1,} & \mbox{ if } x \geq 1,
\end{array}
\right.
\end{equation}
where the regular inverse $F^{-1}$ of $F$ is well defined over
$[F(0), 1]$. $F^*$ is continuous on $(-\infty, \infty)$ and strictly
increasing on $[F(0),1]$.

\noindent\textbf{b.} If the limit $F$ is of Type II, then its
generalized inverse $F^*$, defined in \eqref{G-inverse}, is given by
\begin{equation}\label{case-b}
F^*(x)=\left\{
      \begin{array}{ll}
        0, & \hbox{ if } x<1 \\
        1, & \hbox{ if } x\ge 1.
      \end{array}
    \right.
\end{equation}
This is a degenerate distribution at $x=1$, that is, it induces a
probability measure $\nu=\delta_1$, a delta function at $1$. In this
case,  we have from Theorem~\ref{general} that $\mu_n
\rightsquigarrow \mathrm{Unif}[0,2 \pi)\otimes \delta_1$.  This is
equivalent to that the empirical distribution $\hat\mu_n$ for scaled
eigenvalues converges to the uniform distribution over the unit
circle $|z|=1$ in the complex plane; see Theorem~\ref{typeII}.

\noindent\textbf{c.} When $F$ is of Type III, we have
\begin{equation}\label{case-c}
F^*(x)=\left\{
      \begin{array}{ll}
        0, & \hbox{ if } x<0 \\
        1, & \hbox{ if } x\ge 0.
      \end{array}
    \right.
\end{equation}
This defines a degenerate probability measure $ \delta_0$.  Then the
limit of $\hat\mu_n$ is degenerate at origin in the complex plane.
\end{rem}


The most interesting case to us is the distribution of Type I; see
Theorem~\ref{typeI}. In this case, the normalization constant
$\gamma_n$ should be of the same order as
$\sum^{m_n}_{j=1}\frac{n}{n_j}$, precisely, condition \eqref{vip}
must be true. One can simply take
$\gamma_n=\sum^{m_n}_{j=1}\frac{n}{n_j}$ and calculate the limits
when they exist
\[
c_k:=\lim_{n\to\infty}\frac{\sum^{m_n}_{j=1}\big(\frac{n}{n_j}\big)^k}{\sum^{m_n}_{j=1}\frac{n}{n_j}},
~~~\mbox{ for }k\ge 2.
\]
Then we obtain the limiting distribution $F$ via formula
\eqref{F-rep} with $c_1=1$. Type II and Type III limiting
distributions can be trivially obtained by changing order of
$\gamma_n$.


From Theorem~\ref{typeI},  the limiting empirical distribution of
$\hat\mu_n$ has a support on $F(0)\le |z|\le 1$.  When $F(0)=0$, the
support is the unit disk. When $F(0)>0$,  $\{z: F(0)\le |z|\le 1\}$
is a ring. Since $F(x)$ is right continuous at $x=0$, we have from
\eqref{F-rep} that $F(0)>0$ if and only if $\sum^\infty_{k=1}
\frac{c_k}{k}<\infty$.  Some specific distributions on rings are
given in Examples~\ref{exm2} and \ref{exm3}.

It is interesting to discuss when the empirical distribution
$\mu_n^*$ for linearly scaled eigenvalues converges.  This is
equivalent to the convergence of $\mu_n$ or $\hat\mu_n$ when
$\gamma_n$ is set to be $2$.

When $m_n$ is actually a fixed integer,  Zeng~\cite{Zeng2017}
obtained the liming distribution of $\mu_n^*$ by assuming that
\begin{equation}\label{fixed-m}
\lim_{n\to\infty}\frac{n}{n_j}=:\alpha_j\in [0,1],~~~~2\le j\le m;
 \end{equation}
see Theorem~1.1 in Zeng~\cite{Zeng2017}. By selecting $\gamma_n=2$,
we can verify \eqref{vip} holds, and
\begin{equation}\label{fixed}
    F(x)=x^{1/2}\prod_{j=2}^{m}\big(1-\alpha_j(1-x)\big)^{1/2},~~~x\in (0,1].
\end{equation}
Since $F(0)=0$, the support of the liming distribution of $\mu_n^*$
is always the unit disk $\{z: |z|\le 1\}$.  With additional
constraint $n=n_1\le n_2\le \cdots\le n_m$, it is possible to show
that \eqref{fixed-m} is also necessary for the convergence of $F_n$.


Consider the case $\lim_{n\to\infty}m_n=\infty$. By selecting
$\gamma_n=2$, \eqref{vip} gives the necessary and sufficient
conditions for convergence of $\mu_n^*$.  Again, in this case,
$F(0)=0$ for any limit $F$.  If
$\sum^{m_n}_{r=1}\frac{n}{n_j}\to\infty$, we can only consider the
convergence of the empirical spectral distribution $\mu_n$ for
nonlinearly scaled eigenvalues.

We offer one more comment as a remark before we give some
illustrative examples.

\begin{rem}
 In Theorem~\ref{general}, we have taken $h_n(r)={r}^{2/\gamma_n}/a_n$ for
$r>0$ to re-scale the eigenvalues, where $a_n$ is defined as
$\prod^{m_n}_{j=1}n_j^{1/\gamma_n}$.
 As a matter of fact, if there exist some
sequences $\gamma_n>0$ and $a_n>0$ such that $\mu_n \rightsquigarrow
\mathrm{Unif}[0,2 \pi)\otimes \nu_1$ as $n\to
    \infty$, with $h_n(r)={r}^{2/\gamma_n}/a_n$ and $\nu_1$ being a non-degenerate probability measure, then we can show that
    \begin{equation}\label{lawsoftypes}
    \frac{\gamma_n}{\sum^{m_n}_{j=1}\frac{n}{n_j}}\to
    c~~\mbox{ and }~~\ln  a_n-\frac{\sum^{m_n}_{j=1}\ln  n_j}{\gamma_n}\to d
    \end{equation}
by using the laws of types, where $c\in (0,\infty)$ and $d\in
(-\infty, \infty)$. This implies that there exist some sequences
$\gamma_n>0$ and $a_n>0$ such that $\mu_n \rightsquigarrow
\mathrm{Unif}[0,2 \pi)\otimes \nu_1$ with
$h_n(r)={r}^{2/\gamma_n}/a_n$ and $\nu_1$ being a non-degenerate
probability measure, if and only if $\mu_n \rightsquigarrow
\mathrm{Unif}[0,2 \pi)\otimes \nu$,
 where $\nu$ is a non-degenerate
probability measure with
$h_n(r)=\big(r^2/\prod^{m_n}_{j=1}n_j\big)^{1/\sum^{m_n}_{j=1}n/n_j}$.
Further, the relationship between $\nu$ and $\nu_1$ under condition
\eqref{lawsoftypes} is $\nu_1(-\infty, r]=\nu(-\infty, e^{cd}r^{c}]$
for all $r>0$.

\end{rem}

\begin{exm}\label{exm1}
When these rectangular matrices are actually square matrices, that
is,  $n_1=\cdots=n_{m_n+1}=n$, where $m_n$ is any sequence of
positive integers. Set $\gamma_n=m_n$.  Then \eqref{vip} holds
trivially with $c_k=1$ for all $k\ge 1$.  We have $F(x)=x$, $x\in
(0,1]$. Then $G^*(x)=x$ for $x\in [0,1]$ is the cumulative
distribution function for uniform distribution over $[0,1]$.  This
leads to Theorem~2 in Jiang and Qi~\cite{JiangQi2019}.
\end{exm}

\begin{exm}\label{exm2}
     Let $\{m_n\}$ be positive integers such that
$\lim_{n\to\infty}m_n=\infty$.  Define $n_1=n_{m_n+1}=n$ and assume
$n_2=\cdots=n_{m_n}\sim n\alpha_n$ as $n\to\infty$, where
$\alpha_n\ge 1$. Then
\begin{equation}\label{lambda-app}
\lambda_k(n)=\sum^{m_n}_{r=1}\Big(\frac{n}{n_j}\Big)^k=1+
\frac{m_n}{\alpha_n^k}(1+o(1))  ~~~~\mbox{ as }n\to\infty
\end{equation}
for $k\ge 1$.

Assume $\lim_{n\to\infty}\alpha_n=:\alpha\in [1,\infty)$. By taking
$\gamma_n=2m_n$,  we see that \eqref{vip} holds with $c_1=\frac12$,
and $c_1=\frac12\alpha^{-k}$ for $k\ge 1$. We have
\[
F(x)=\Big(1-\frac{1}{\alpha}(1-x)\Big)^{1/2}~~~~~~x\in (0,1].
\]
Then we obtain
\[
 F^{*}(x)=
\left\{
\begin{array}{ll}
{0,} &\mbox{ if } x<\beta; \\
\frac{x^2-\beta^2}{1-\beta^2}, & \mbox{ if } x \in [\beta,1); \\
{1,} & \mbox{ if } x \geq 1
\end{array}
\right.
\]
with $\beta=(1-\frac{1}{\alpha})^{1/2}$. The density function of
$F^*$ is given by $f^*(x)=\frac{2x}{1-\beta^2}I(\beta\le x\le 1).$
According to Theorem~\ref{typeI}, $\hat\mu_n\rightsquigarrow
\mathrm{Unif}(\beta\le |z|\le 1)$. The limit is a uniform
distribution on the ring $\beta\le |z|\le 1$ if $\beta\in (0,1)$,
and a uniform distribution on the unit disk if $\beta=0$.
\end{exm}

\begin{exm}\label{exm3} In Example~\ref{exm2}, we assume
$\lim_{n\to\infty}\alpha_n=\infty$.

\noindent\textbf{(a)}. Consider the case
$\lim_{n\to\infty}\frac{m_n}{\alpha_n}=\infty$. With selecting
$\gamma_n=\sum^{m_n}_{r=1}\frac{n}{n_j}$, we have $c_1=1$ and
$c_k=0$ for all $k\ge 2$. Then we have $F(x)=\exp(x-1)$ for $x\in
(0,1]$, yielding
\[
    F^{*}(x)=
\left\{
\begin{array}{ll}
{0,} &\mbox{ if } x<e^{-1}; \\
1+\ln x, & \mbox{ if } x \in [e^{-1},1); \\
{1,} & \mbox{ if } x \geq 1.
\end{array}
\right.
\]
It follows from Theorem~\ref{typeI} that $\hat\mu_n\rightsquigarrow
\mu$, where $\mu$ has a density function $\frac{1}{2\pi
|z|^2}I(e^{-1}\le |z|\le 1)$.

\noindent\textbf{(b)}.  Consider the case
$\lim_{n\to\infty}\frac{m_n}{\alpha_n}=\gamma\in [0,\infty)$.  It
follows from \eqref{lambda-app} that
$\lim_{n\to\infty}\lambda_1(n)=1+\gamma$, and
$\lim_{n\to\infty}\lambda_k(n)=1$ for $k\ge 2$. This is the case we
can establish the limiting law for $\mu_n^*$, the empirical
distribution for linearly scaled eigenvalues, as defined in
\eqref{mu*}. By selecting $\gamma_n=2$, we have
\[
F(x)=x^{1/2}\exp\big(\frac{\gamma}{2}(x-1)\big), ~~~~~x\in (0,1].
\]
Let $f^*$ denote the density of $F^*=F^{-1}$ on $(0,1))$. We have
$\mu_n^*\rightsquigarrow \mu$, where $\mu$ has a density function
$\frac{f^*(|z|)}{2\pi |z|}$.
\end{exm}

\begin{figure}[!h]
\centering \caption{Scatter plots for product matrices:
$n_1=n_{m+1}=n$, $n_2=\cdots=n_m=2n$, $\gamma_n=2m$}
\includegraphics[height=3.0in,width=5in]{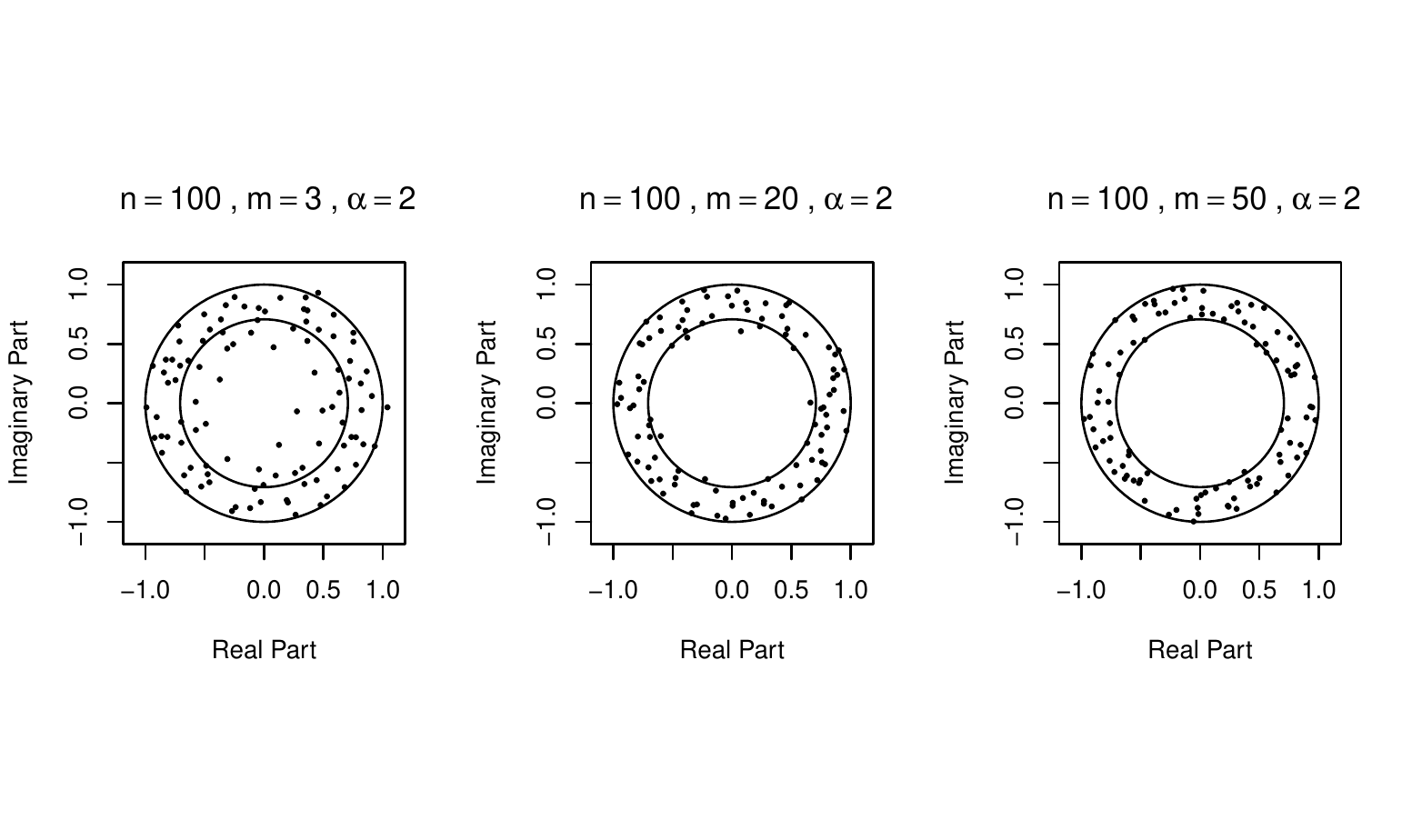}
\label{scatter1}
\end{figure}
\begin{figure}[!h]
\centering \caption{Scatter plots for product matrices:
$n_1=n_{m+1}=n$, $n_2=\cdots=n_m=2n$, $\gamma_n=2m$}
\includegraphics[height=3.0in,width=5in]{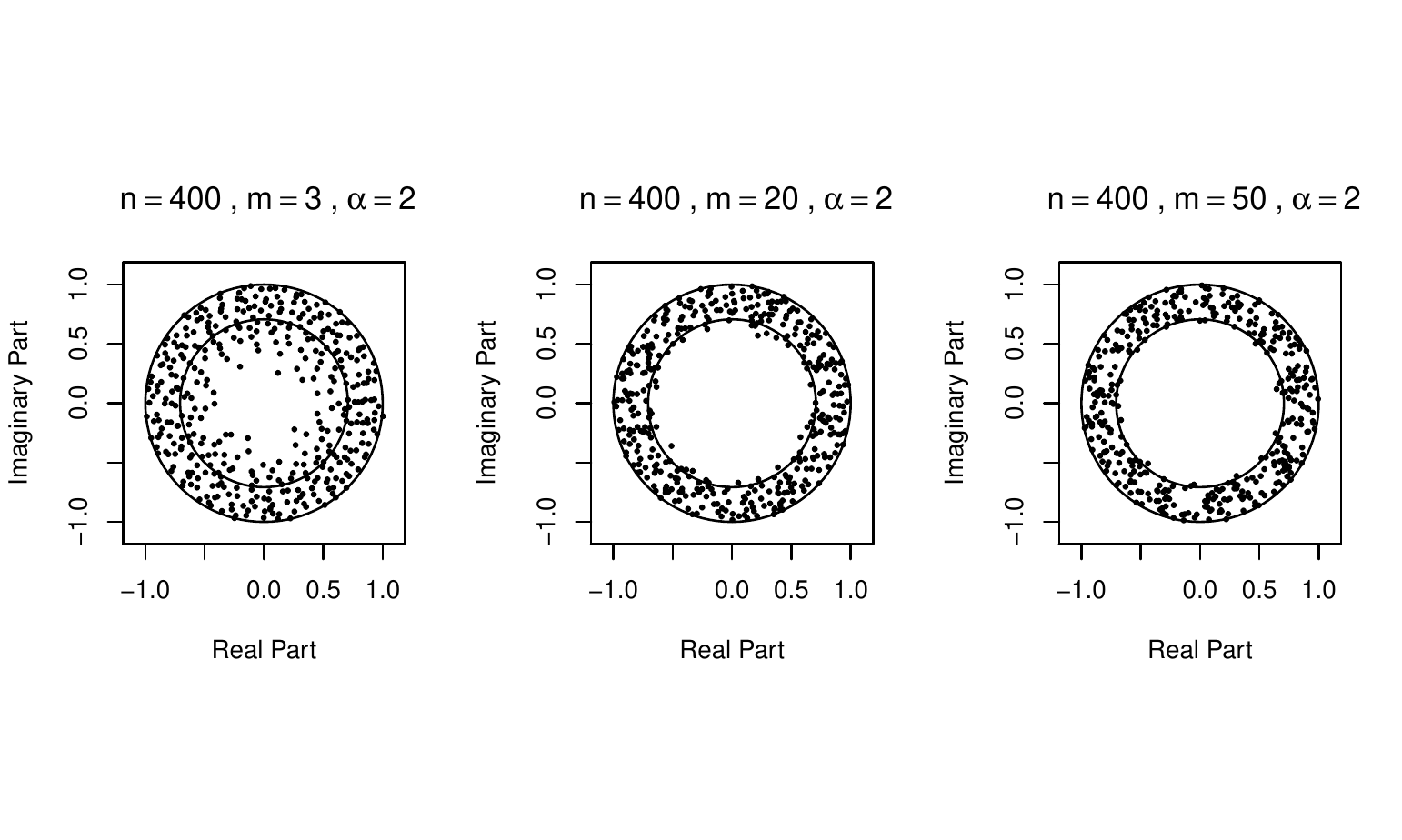}
\label{scatter2}
\end{figure}

To conclude this section, we carry out a simulation study by using
the setup in Example~\ref{exm2}. We select $\alpha_n=\alpha=2$,
$n_2=\cdots=n_m=2n$,and $\gamma_n=2m$.   Theoretically, if $m$ is
large, the empirical spectral distribution for the nonlinearly
scaled eigenvalues is approximately uniformly distributed on the
ring $\{\frac{\sqrt{2}}2\le |z|\le 1\}$.  For each of $n=100$ and
$n=400$, we select $m=3$, $m=20$ and $m=50$ in order to see how well
these scaled eigenvalues fit into the ring with the change in the
value of $m$.  The scatter plots for the scaled eigenvalues when
$n=100$ and $n=400$ are given in Figures~\ref{scatter1} and
~\ref{scatter2}, respectively. From the two figures, we see that
most of the scaled eigenvalues are already falling within the ring
$\{\frac{\sqrt{2}}2\le |z|\le 1\}$ when $m=20$.

\section{Proofs}\label{proofs}

The lemmas~\ref{nonlinear} and ~\ref{lemJQ} below play a very
important role in the proofs of our main results.

\begin{lemma}\label{nonlinear} (Theorem 1 in Jiang and Qi~\cite{JiangQi2019}). Let $\varphi(x)\geq 0$ be a measurable  function defined on $[0, \infty).$ Assume the density  of $(Z_1, \cdots, Z_n)\in \mathbb{C}^n$ is proportional to  $\prod_{1\leq j < k
\leq n}|z_j-z_k|^2\cdot \prod_{j=1}^n\varphi(|z_j|)$. Let $Y_1,
\cdots, Y_n$ be independent r.v.'s such that the density of $Y_j$ is
proportional to  $y^{2j-1}\varphi(y)I(y\geq 0)$ for every $1\leq
j\leq n.$ Let $\mu_n$,$\nu_n$ and $\nu_n^*$ be defined as in
\eqref{mun} and \eqref{nun*}, respectively. If $\{h_n\}$ are
measurable functions such that $\nu_n\rightsquigarrow\nu$ for some
probability measure $\nu$, then $\mu_n \rightsquigarrow \mu$ with
$\mu=\mathrm{Unif}[0, 2\pi]\otimes\nu$ .
Taking $h_n(r)=r/a_n$, the conclusion still holds if ``$(\mu_n,
\nu_n, \mu, \nu)$" is replaced by  ``$(\mu_n^*, \nu_n^*, \mu^*,
\nu^*)$" where $\mu^*$ is the distribution of $Re^{i\Theta}$ with
$(\Theta, R)$ having the law of  $\mathrm{Unif}[0,
2\pi]\otimes\nu^*$.
\end{lemma}

Let $Y_1, \cdots, Y_n$ be the independent random variables
determined in Lemma~\ref{nonlinear} under model \eqref{model}. Let
$\{s_{j,r}, ~1 \leq j \leq n,~ 1 \leq r \leq m_n\}$ be independent
random variables and $s_{j,r}$ follow a Gamma($l_r+j$) with density
function $y^{l_{r}+j-1} e^{-y} I_{y>0} /
\Gamma\left(l_{r}+j\right)$. Set
\begin{equation}\label{Tj}
 T_{j}=\prod_{r=1}^{m_n} s_{j, r}, ~~ 1 \leq j \leq n.
\end{equation}

\begin{lemma}\label{lemJQ} (Lemma 4 in Jiang and Qi~\cite{JiangQi2019})
Suppose $\{h_n(x);\, n\geq 1\}$ are measurable functions defined on
$[0,\infty)$ and $\nu_n$'s are defined as in \eqref{mun}. Let $Y_1,
\cdots, Y_n$ be as in Lemma~\ref{nonlinear} and $\nu$ be a
probability measure on $\mathbb{R}.$ Then $\nu_n\rightsquigarrow\nu$
if and only if
 \[
\lim_{n\to\infty}\frac{1}{n}\sum^n_{j=1}\mathbb{P}(h_n(Y_j)\le
r)=H(r)
\]
for every continuity point $r$ of $H(r)$, where $H(r):=\nu((-\infty,
r]),\, r \in \mathbb{R}$.
\end{lemma}

The results in the following lemma are summarized  from Lemmas~2.2
and 2.3 from  Zeng~\cite{Zeng2017}.

\begin{lemma}\label{zeng} (Zeng~\cite{Zeng2017}) We have
\begin{equation}\label{eq1}
\mathbb{P}(T_{1} \leq x) \geq \mathbb{P}(T_{2} \leq x) \geq \cdots
\geq \mathbb{P}(T_{n} \leq x)
\end{equation}
for any $x \in[0, \infty)$, $(Y_1^2, \cdots,
Y_n^2)\overset{d}=(T_1,\cdots, T_n)$, and
\begin{equation}\label{eq2}
g(T_{1}, \cdots, T_{n})\overset{d}=g(|z_{1}|^{2}, \cdots,
|z_{n}|^{2})
\end{equation}
for any symmetric function $g(t_1,\cdots,t_n)$, where $\overset{d}=$
denotes equality in distribution.
\end{lemma}


Before we prove Theorem~\ref{type}, we need to introduce more
notation and preliminary results.

Define
\begin{equation}\label{lambda}
\lambda_k(n)=\sum^{m_n}_{j=1}(\frac{n}{n_j})^k,  ~~~k\ge 1
\end{equation}
and
\begin{equation}\label{theta}
\theta_k(n)=\frac{\lambda_k(n)}{\lambda_1(n)}, ~~k\ge 1.
\end{equation}
Note that $\lambda_k(n)\ge 1$ since $n_1=n$. Since $n\le n_j$ for
all $1\le j\le m_n$ we have $\lambda_k(n)$ is non-increasing in
$k\ge 1$ for each $n$. Thus, we have $\theta_k(n)$ is non-increasing
in $k\ge 1$, implying that $0<\theta_k(n)\le \theta_1(n)=1$ for
$k\ge 1$ and $n$.

We define a sequence of new distribution functions as follows
\begin{equation}\label{Gn}
    G_{n}(x)=
    \Big(\prod_{j=1}^{m_n} (1-\frac{n}{n_j}(1-x))\Big)^{1 /\lambda_1(n)},~~~~  x \in[0,1]
\end{equation}
These distributions are obtained by letting $\gamma_n=\lambda_1(n)$
in equation \eqref{Fn}.  We have
\begin{equation}\label{FG}
F_n(x)=G_n^{\lambda_1(n)/\gamma_n}(x), ~~~x\in [0,1].
\end{equation}

Set $g_n(x)=\ln(G_n(x))$, $x\in(0,1]$. Then $g_n(x)\le 0$ for $x\in
(0,1]$. Using Taylor's expansion $\ln (1-t)=-\sum^\infty_{k=
1}\frac{t^k}{k}$, $|t|<1$, we have for $x\in(0,1]$
    \begin{equation}\label{gn}
    g_n(x)=-\sum_{k=1}^\infty \frac{\theta_k{(n)}}{k}(1-x)^k.
    \end{equation}
Note that $\theta_1(n)=1$ and $0<\theta_k(n)\le 1$. We have the
following inequalities
    \begin{equation}\label{inequalities}
 1-x\leq-g_{n}(x)\leq \sum_{k=1}^\infty \frac{1}{k}(1-x)^k =-\ln
(x),~~~x\in(0,1].
    \end{equation}

The probability distribution $G_n$ over $[0,1]$ is defined via
\eqref{Gn}, and the same expression can not be extended beyond the
interval $[0,1]$. The function $g_n(x)$, as the logarithm of $G_n$,
has expansion \eqref{gn} over $(0, 1]$ only.  However, $g_n$ can be
extended to a region in the complex plane via the expression on the
right-hand side of \eqref{gn}. Now we fix $0<\delta<1$. For any
complex number $z$ such that $|z-1|\leq \delta$, we have from
\eqref{inequalities} that
\begin{equation}\label{unif-bound}
  \sum_{k=1}^\infty \left|-\frac{\theta_k{(n)}}{k}(1-z)^k\right|\leq \sum_{k=1}^\infty \frac{\theta_k{(n)}}{k}\delta^k=|g_n(1-\delta)|\le
\ln(\frac{1}{1-\delta}).
\end{equation}
Therefore, we can extent $g_n(x)$ to be a complex analytic function
on disk $D=\{z\in\mathbb{C}: |z-1|<1\}$, namely
    \begin{equation}\label{equ:taylor g_n(x)}
    g_n(z)=-\sum_{k=1}^\infty \frac{\theta_k{(n)}}{k}(1-z)^k.
    \end{equation}


\begin{lemma}\label{complexlimit} (Theorem 10.28 in Rudin~\cite{RudinComplex}) Suppose $f_j$ is analytic on
open set $\Omega\subset \mathbb{C}$ for $j=1,2,\cdots,$ and
$f_j\rightarrow f$ uniformly on each compact subset of $\Omega$.
Then $f$ is analytic on $\Omega$, and $f'_j\rightarrow f'$ uniformly
on any compact subset on $\Omega$.
\end{lemma}

\begin{lemma}\label{sufficiency}
Assume $\{n_s\}$ is a subsequence of $\{n\}$ such that
$\lim_{s\to\infty}\theta_k(n_s)=a_k\in[0,1]$ for all $k\ge 2$. Set
$a_1=1$. Then
\begin{equation}\label{Gn-limit}
\lim_{s\to\infty}G_{n_s}(x)=G(x), ~~~x\in (0,1],
\end{equation}
where $G$ is a distribution function given by
\begin{equation}\label{Grep}
G(x)=\exp(-\sum^\infty_{k=1}\frac{a_k}{k}(1-x)^k), ~~~x\in (0,1].
\end{equation}
$0<G(x)<1$ is analytic and strictly increasing over $(0,1)$.
\end{lemma}

\noindent{\it Proof.}  In the proof  we will use index $n$ instead
of $n_s$ for the sake of brevity.

For each $\delta\in(0,1)$, set $K_\delta:=\{z\in \mathbb{C}:
 |z-1|\leq\delta\}$. It follows form \eqref{unif-bound} that
 $g_n(z)$ is uniformly bounded on $K_\delta$.

Set $g(z)=-\sum_{k=1}^\infty \frac{a_k}{k}(1-z)^k$. The radius of
convergence of $g(z)$ satisfies
$$
R=\frac{1}{\limsup_{n\rightarrow \infty} \left(
\frac{a_k}{k}\right)^{1/k} }\geq\frac{1}{\limsup_{n\rightarrow
\infty} \left( \frac{1}{k}\right)^{1/k} } =1,
$$
i.e. $g(z)$ is well defined on disk $D=\{z\in\mathbb{C}: |z-1|<1\}$.
For each $\delta\in(0,1)$, we have
\begin{align*}
  \sup_{z\in K_\delta} |g_n(z)-g(z)|\leq \sum_{k=1}^N\frac{|a_k{(n)}-a_k|}{k}\delta^k
+\sum_{k=N+1}^\infty \frac{\delta^k}{k}\leq
\sum_{k=1}^N|a_k{(n)}-a_k|+\frac{\delta^{N+1}}{1-\delta},
\end{align*}
which implies
\begin{eqnarray*}
\limsup_{n\rightarrow \infty} \sup_{z\in K_\delta}
|g_n(z)-g(z)|&\le& \limsup_{N\rightarrow
\infty}\limsup_{n\rightarrow \infty}
\Big(\sum_{k=1}^N|a_k{(n)}-a_k|+\frac{\delta^{N+1}}{1-\delta}\Big)\\
&\leq &  \limsup_{N\rightarrow
\infty}\frac{\delta^{N+1}}{1-\delta}\\
&=&0,
\end{eqnarray*}
that is, $g_n(z)$ converges to $g(z)$ uniformly on $K_\delta$.

Since $G_n(x)=\exp(g_n(x))$ for $x\in (0,1]$, we have
\eqref{Gn-limit} with $G(x)=\exp(g(x))$ for $x\in (0,1]$. Note that
$g(x)$ is analytic,  $g(x)<0$ and is strictly increasing for $x\in
(0,1)$, we have $0<G(x)<1$ is analytic and strictly increasing over
$(0,1)$.
 \hfill$\blacksquare$

\begin{lemma}\label{necessity}
Assume $\{n_s\}$ is a subsequence of $\{n\}$ such that  $G_{n_s}$
converges weakly to a distribution $G$, then
$\lim_{s\to\infty}\theta_k(n_s)=:a_k\in[0,1]$ for all $k\ge 2$, and
$G$ has a representation \eqref{Grep}.
\end{lemma}

\noindent{\it Proof.} Note that $0\leq \theta_k(n)\leq\theta_1(n)=1$
for all $n\ge 1$ and $k\geq 2$. By the diagonal argument, for every
subsequence of $\{n\}$, we can find its further subsequence along
which $\theta_k(n)$ has a subsequential limit in $[0,1]$ for all
$k\ge 2$.

We aim to show that $\lim_{s\to\infty}\theta_k(n_s)$ exists for all
$k\ge 2$. If the conclusion is not true, then for some $k\ge 2$, say
$k_0$, such that the limit of $\theta_{k_0}(n_s)$ doesn't exist.
Then there exist two subsequences of $\{n_s\}$, say $\{n_{s'}\}$ and
$\{n_{s''}\}$, such that
\begin{equation}\label{different-lim}
\lim_{s'\to\infty}\theta_{k_0}(n_{s'})=a \ne
b=\lim_{s''\to\infty}\theta_{k_0}(n_{s''}).
\end{equation}
By the diagonal argument, we can find a further subsequence of
$\{n_{s'}\}$, along which $\theta_k(n_{s'})$ has a subsequential
limit $a_k\in [0,1]$ for each $k\ge 2$ with $a_{k_0}=a$. By
Lemma~\ref{sufficiency} we have
\begin{equation}\label{G1}
G(x)=\exp(-\sum^\infty_{k=1}\frac{a_k}{k}(1-x)^k), ~~~x\in (0,1].
\end{equation}
since any subsequential limit of $G_{n_s}(x)$ is equal to $G(x)$ in
$x\in (0,1]$. Similarly,  we can find a further subsequence of
$\{n_{s''}\}$, along which $\theta_k(n_{s''})$ has a subsequential
limit $b_k\in [0,1]$ for each $\ge k$ with $b_{k_0}=b$. Again, using
Lemma~\ref{sufficiency} we have
\begin{equation}\label{G2}
G(x)=\exp(-\sum^\infty_{k=1}\frac{b_k}{k}(1-x)^k), ~~~x\in (0,1].
\end{equation}
By combining \eqref{G1} and \eqref{G2}, we have
\[
\sum^\infty_{k=1}\frac{a_k}{k}(1-x)^k=\sum^\infty_{k=1}\frac{b_k}{k}(1-x)^k,~~~x\in
(0,1].
\]
Therefore, we $a_k=b_k$ for all $k\geq 2$, which contradicts
$a_{k_0}=a\ne b= b_{k_0}$ from \eqref{different-lim}. This proves
the lemma. \hfill$\blacksquare$

\vspace{10pt}

\noindent{\it Proof of Theorem~\ref{type}.} First, we assume $F_n$
converges weakly to a distribution function $F$. We will show $F$
must be of one of the three types given in Theorem~\ref{type}.

Review the definitions of $\lambda_k(n)$ and $\theta_k(n)$ in
\eqref{lambda} and \eqref{theta}, respectively.

We consider the sequence $\lambda_1(n)/\gamma_n$. At this moment, we
don't know yet whether $ r_n:=\lambda_1(n)/\gamma_n$ has a limit. We
assume that $\{n_s\}$ is any subsequence of $\{n\}$ such that
\begin{equation}\label{sublimit}
\lim_{s\to\infty}
r_{n_s}=\lim_{s\to\infty}\frac{\lambda_1(n_s)}{\gamma_{n_s}}=c_1\in
[0,\infty].
\end{equation}
We consider the following three cases individually:  $c_1\in (0,1)$,
$c_1=\infty$, and $c_1=0$.  From \eqref{FG} we have that
\begin{equation}\label{GF}
G_n(x)=F_n^{1/ r_n}(x)~~~~~x\in (0,1].
\end{equation}

\noindent{\it Case 1.} $c_1\in (0,\infty)$ in \eqref{sublimit}.

In this case,  we see that $G_{n_s}$ converges weakly to
$G(x)=F^{1/c_1}(x)$.  By applying Lemma~\ref{necessity}, we have
$\lim_{s\to\infty}\theta_k(n_s)=:a_k\in[0,1]$ for all $k\ge 2$, and
$G$ has a representation \eqref{Grep}, which implies  $F$ has a
representation \eqref{F-rep} with $c_k=c_1a_k$ for $k\ge 2$.  This
shows that $F$ is of type I.

\noindent{\it Case 2.} $c_1=0$ in \eqref{sublimit}.

In view of \eqref{FG}, \eqref{gn} and \eqref{inequalities}, we have
for any $x\in (0,1)$
\[
1\ge F_{n_s}(x)=G_n^{ r_{n_s}}(x)=\exp( r_{n_s}g_n(x))\ge \exp(
r_{n_s}\ln(x))\to 1
\]
as $s\to\infty$, which implies $F(x)=1$ for $x\in (0,1)$ and thus
$F$ is of type II.

\noindent{\it Case 3.} $c_1=\infty$ in \eqref{sublimit}.

Using the same equations as in the proof for {\it Case 2}, we have
for any $x\in (0,1)$
\[
0\le F_{n_s}(x)=G_n^{ r_{n_s}}(x)=\exp( r_{n_s}g_n(x))\le \exp(
r_{n_s}(1-x))\to 0
\]
as $s\to\infty$, which implies $F(x)=0$ for $x\in (0,1)$ and thus
$F$ is of type III.

We have proved that there are only three types of limiting
distributions for $F_n$.   Next, we will show the necessary and
sufficient conditions in parts (a), (b), and (c).

Sufficiency for parts (b) and (c) has been proved. In fact, for part
(b), condition \eqref{zero} must be true when $F$ is of Type II,
otherwise, there exists a subsequential limit $c_1$ of $ r_n$ with
$c_1\in (0,\infty)$ or $c_1=\infty$, such that $F$ is of Type I or
Type III, respectively,  yielding a contradiction.  A similar
argument can be used to show \eqref{infinity} in part (c).

Finally,  we need to prove part (a).  The sufficiency has been
proved in {\it Case 1} above.  Assume $F_n$ converges weakly to $F$,
which is of Type I.  We show \eqref{vip}, or equivalently, we show
the following statements

\noindent{\it Statement 1:} $r_n=\lambda_1(n)/\gamma_n$ has a limit
$c_1\in (0,\infty)$;

\noindent{\it Statement 2:} For any $k\ge 2$, $\theta_{k}(n)$ has a
finite limit.

If {\it Statement 1} is not true, then there are subsequences of
$\{n\}$, say, $\{n_{s}\}$ and $\{n_{s'}\}$ such that
\[
\lim_{s\to\infty}r_{s}=a\ne b=\lim_{s'\to\infty}r_{s'},
\]
and $a, b\in (0,\infty)$.  Any subsequential limit of $r_n$ must be
a finite positive number since $F$ is of Type I.

From \eqref{GF}, we have $G_{n_s}$ converges weakly to $F^{1/a}$ and
$G_{n_{s'}}$ converges weakly to $F^{1/b}$.  Then it follows from
Lemma~\ref{necessity} that
\[
F^{1/a}(x)=\exp(-\sum^\infty_{k=1}\frac{a_k}{k}(1-x)^k), ~~~x\in
(0,1]
\]
and
\[
F^{1/b}(x)=\exp(-\sum^\infty_{k=1}\frac{b_k}{k}(1-x)^k), ~~~x\in
(0,1]
\]
where $a_1=b_1=1$ and $a_k, b_k\in [0,1]$ for all $k\ge 2$.   We
conclude that
\[
a\sum^\infty_{k=1}\frac{a_k}{k}(1-x)^k=b\sum^\infty_{k=1}\frac{b_k}{k}(1-x)^k,~~~x\in
(0,1].
\]
Since the functions on both sides of the above equation are
analytic, their first derivatives at $x=1$ must be the same, which
leads to $a=b$,  contradictory to the assumption $a\ne b$.
Therefore, {\it Statement 1} is true, that is, $r_n$ has a limit in
$(0,\infty)$.

Given $\lim_{n\to\infty}r_n=c_1$,  from \eqref{GF} we have $G_{n}$
converges weakly to $F^{1/c_1}$.  Again, by using
Lemma~\ref{necessity}, we have $a_k:=\lim_{n\to\infty}\theta_k(n)$
exists for all $k\ge 2$ and $a_k\in [0,1]$.   This proves {\it
Statement 2}.  The proof of the theorem is completed.
\hfill$\blacksquare$

 \vspace{10pt}

The following result is an extension of Lemma 2.3 in Zeng
~\cite{Zeng2017}. We allow $m_n$ to change with $n$.

\begin{lemma}\label{Tnlimit}
    Assume $\{m_n\}$ is a sequence of positive integers. Then
    \begin{equation}\label{limit1}
    \frac{1}{\lambda_1(n)} \ln \frac{T_{[\mathrm{nx}]}}{\prod_{r=1}^{m_n}(l_r+n)}-\ln G_n(x) \stackrel{p}{\rightarrow} 0,~~ x \in
(0,1],
    \end{equation}
where $[nx]$ denotes the integer part of $nx$,  $T_j$ is defined in
\eqref{Tj}, $\lambda_1(n)$ is defined in \eqref{lambda}, and $G_n$
is defined in \eqref{Gn}.
\end{lemma}

\noindent{\it Proof.} We have $\ln T_{j}=\sum_{r=1}^{m} \ln s_{j,
r}$ for $j \geq 1$. Since $s_{j,r}$ has a Gamma($l_r+j$)
distribution, we have
    \[
   \mu_{j, r}=\mathbb{E}\left(s_{j, r}\right)=l_{r}+j,~~{
\Var}\left(s_{j, r}\right)=l_{r}+j
    \]
and the moment generating function of $\ln s_{j, r}$ is
    \[
        m_{j}(t)=\mathbb{E}\big(e^{t \ln s_{j, r}}\big)=\frac{\Gamma\left(l_{r}+j+t\right)}{\Gamma\big(l_{r}+j\big)},~~~ t>-(l_r+j).
    \]
    It follows that
    \[
        \mathbb{E}\big(\ln s_{j, r}\big)=\frac{d}{d t} m_{j}(t)\Big|_{t=0}=\frac{\Gamma^{\prime}\left(l_{r}+j\right)}{\Gamma\left(l_{r}+j\right)}=\psi\left(l_{r}+j\right),
    \]
    where $\psi(x)=\Gamma^{\prime}(x) / \Gamma(x)$ is a digamma function.
    Thus, we have
    \begin{equation}\label{ElogT}
        \mathbb{E}\big(\ln T_{j}\big)=\sum_{r=1}^{m_n} \mathbb{E}\left(\ln s_{j, r}\right)=\sum_{r=1}^{m_n}
\psi\left(l_{r}+j\right).
    \end{equation}

Set $\eta(t)=t-1-\ln t$ for $t>0$.  Then $\eta(t)\ge 0$ for $t>0$.

Trivially, we have for any $1\le j\le n$
    \begin{equation}\label{logTj}
   \ln \frac{T_{j}}{\prod_{r=1}^{m_n}\mu_{j,r}}-\sum_{r=1}^{m_n}
\ln
(\frac{l_r+j}{l_r+n})=\sum_{r=1}^{m_n}(\frac{s_{j,r}}{\mu_{j,r}}-1)-\sum_{r=1}^{m_{n}}\eta(\frac{s_{j,r}}{\mu_{j,r}}).
    \end{equation}


It can be seen that for $1\le j\le n$
    \[
        {\Var}\left(\sum_{r=1}^{m_n}\Big(\frac{s_{j, r}}{\mu_{j, r}}-1\Big)\right)
        =\sum_{r=1}^{m_n} \frac{{\Var}\left(s_{j, r}\right)}{\mu_{j,
r}^{2}}=\sum_{r=1}^{m_n}\frac{1}{l_{r}+j}\le\sum_{r=1}^{m_n}\frac{n}{j(l_r+n)}=\frac{\lambda_1(n)}{j}.
    \]

Fix $x\in (0,1]$.  Set $j=j_n=[nx]$. Then
 \[
    \Var\Big(\frac{1}{\lambda_1(n)}\sum_{r=1}^{m_n}(\frac{s_{j_n, r}}{\mu_{j_n,
r}}-1)\Big)\le \frac{1}{\lambda_1(n)j_n}\le\frac{1}{j_n} \to 0
 \]
as $n\to\infty$. By Chebyshev inequality, we obtain
    \begin{equation}\label{term1}
    \frac{1}{\lambda_1(n)}\sum_{r=1}^{m_n}(\frac{s_{j_n, r}}{\mu_{j_n,
r}}-1) \stackrel{p}{\rightarrow} 0.
    \end{equation}

From \eqref{logTj} and \eqref{ElogT} we have
    \[
    \mathbb{E}\Big(\sum_{r=1}^{m_n} \eta\big(\frac{s_{j, r}}{\mu_{j, r}}\big)\Big)
            =\sum_{r=1}^{m_n} \ln \mu_{j, r}-\mathbb{E} \sum_{r=1}^{m_n} \ln s_{j, r}
           =-\sum_{r=1}^{m_n}\big(\psi(l_{r}+j)-\ln (l_{r}+j)\big).
     \]
We need the following approximation for $\psi$
\begin{equation*}
 \psi(t)-\ln t=-\frac{1}{2 t}+O\Big(\frac{1}{t^{2}}\Big) \mbox{ as }t\to
\infty;
\end{equation*}
See, e.g., Formula 6.3.18 in Abramowitz and
Stegun~\cite{Abramowitz1972}. With $j_n=[nx]$, we have
\[
\frac{1}{\lambda_1(n)}\mathbb{E}\Big(\sum_{r=1}^{m_n}
\eta\big(\frac{s_{j, r}}{\mu_{j, r}}\big)\Big)
=\frac{O(1)}{\lambda_1(n)}\sum_{r=1}^{m_n}\frac{1}{l_{r}+j_n}\le
\frac{O(1)}{\lambda_1(n)}\frac{\lambda_1(n)}{j_n}=O(\frac1{j_n})\to
0
 \]
as $n\to\infty$, which implies
    \begin{equation}\label{term2}
    \frac{1}{\lambda_1(n)}\sum_{r=1}^{m_n} \eta\left(\frac{s_{j, r}}{\mu_{j, r}}\right) \stackrel{p}{\rightarrow} 0
    \end{equation}
by Chebyshev inequality since $\sum_{r=1}^{m} \eta\left(\frac{s_{j,
r}}{\mu_{j, r}}\right)\ge 0$.   Therefore, combining (\ref{logTj}),
(\ref{term1}) and (\ref{term2}), we obtain
    \begin{equation}\label{limit2}
     \frac{1}{\lambda_1(n)} \ln \frac{T_{[\mathrm{nx}]}}{\prod_{r=1}^{m_n}(l_r+[nx])}-\ln G_n(\frac{[nx]}{n}) \stackrel{p}{\rightarrow} 0,
 ~~~ x \in (0,1].
    \end{equation}

From \eqref{gn}, we have
$0<g_n'(t)=\sum^\infty_{r=1}\theta_k(n)(1-t)^{k-1}\le \frac1{t}$ for
$0<t\le 1$, and hence,
\[
|\ln G_n(\frac{[nx]}{n})-\ln G_n(x)|=|g_n(\frac{[nx]}{n})-g_n(x)|\le
\sup_{\frac{[nx]}{n}\le t\le x}g'_n(t)|\frac{[nx]}{n}-x|\le
\frac{1}{[nx]}\to 0
\]
as $n\to\infty$. This, coupled with \eqref{limit2}, yields
\eqref{limit1}.
     \hfill$\blacksquare$

\vspace{10pt}

\noindent{\it Proof of Theorem~\ref{general}.}   Assume $F_n$
converges weakly to a distribution $F$. The conclusion in the
theorem follows from Lemma~\ref{nonlinear} and Lemma~\ref{lemJQ} if
we can prove
\[
\lim_{n\to\infty}\frac{1}{n}\sum^n_{j=1}\mathbb{P}(h_n(Y_j)\le
y)=F^*(y)
\]
for every continuity point $y$ of $F^*$. According to \eqref{eq2},
it is equivalent to show
\begin{equation}\label{expectation1}
\lim_{n\to\infty}\frac{1}{n}\sum^n_{j=1}\mathbb{P}\Big(\frac{1}{a_n}T_j^{1/\gamma_n}\le
y\Big)=F^*(y)
\end{equation}
for every continuity point $y$ of $F^*$.   Since
$a_n=\prod^{m_n}_{r=1}n_r^{1/\gamma_n}=\prod^{m_n}_{r=1}(l_r+n)^{1/\gamma_n}$,
we have
\begin{equation}\label{eq30}
\mathbb{P}\Big(\frac{1}{a_n}T_j^{1 / \gamma_{n}} \leq y\Big)=
\mathbb{P}\Big(\frac{1}{\gamma_{n}} \ln
\frac{T_{j}}{\prod_{r=1}^{m_{n}}\left(l_{r}+j\right)} \leq \ln
y\Big),~~~y>0.
\end{equation}
We also have the following two inequalities
\begin{equation}\label{estimates}
\frac{j}n\mathbb{P}\Big(\frac{1}{a_n}T_j^{1/ \gamma_{n}} \leq
y\Big)\le \frac{1}{n} \sum_{j=1}^{n}
\mathbb{P}\Big(\frac{1}{a_n}T_j^{1/ \gamma_{n}} \leq y\Big)\le
\frac{j}{n}+\mathbb{P}\Big(\frac{1}{a_n}T_j^{1/ \gamma_{n}} \leq
y\Big),
\end{equation}
which follow from the monotonicity in \eqref{eq1} directly.

\noindent{\underline{\it Case 1}.} Assume $F$ is of Type I.

In this case, $F$ is strictly increasing in $[0,1]$ with $F(x)>0$
for any $x\in (0,1]$ and $\lim_{n\to\infty}F_n(x)=F(x)$, and $F^*$
is given by \eqref{case-a}. We note that $\lambda_1(n)/\gamma_n$
converges to a non-zero constant from Theorem~\ref{type}. From
\eqref{limit1}, we have for any $x\in (0,1]$
\[ \frac{1}{\gamma_n}\ln
\frac{T_{[\mathrm{nx}]}}{\prod_{r=1}^{m_n}(l_r+n)}-\ln F_n(x)=
\frac{\lambda_1(n)}{\gamma_n}\Big(\frac{1}{\lambda_1(n)} \ln
\frac{T_{[\mathrm{nx}]}}{\prod_{r=1}^{m_n}(l_r+n)}-\ln G_n(x))\Big)
\stackrel{p}{\rightarrow} 0,
 \]
yielding
  \begin{equation}\label{limit3}
\frac{1}{\gamma_n}\ln
\frac{T_{[\mathrm{nx}]}}{\prod_{r=1}^{m_n}(l_r+n)}-\ln F(x)
\stackrel{p}{\rightarrow} 0.
    \end{equation}

From \eqref{eq30} we get
\begin{equation}\label{eq3}
\mathbb{P}\Big(\frac{1}{a_n}T_j^{1 / \gamma_{n}} \leq y\Big)=
\mathbb{P}\Big(\frac{1}{\gamma_{n}} \ln
\frac{T_{j}}{\prod_{r=1}^{m_{n}}\left(l_{r}+j\right)}-\ln F(x) \leq
\ln \frac{y}{F(x)}\Big)
\end{equation}
for $y>0$ and any $x$ with $F(x)>0$.

Now we are ready to show \eqref{expectation1} when $F(0)<y<1$,
$y\leq F(0)$,  and $y\geq 1$.

We first assume $F(0)<y<1$.  Let $\delta \in(0,1)$ be any given
number such that
\[F(0)< y-\delta<y<y+\delta<1.\]
Then $0<F^{*}(y-\delta)<F^{*}(y)<F^{*}(y+\delta)<1$.

By setting $x=F^{*}(y+\delta)$, $j=[nx]$ in \eqref{estimates} and
\eqref{eq3} and using \eqref{limit3},  we have
\begin{eqnarray*}
&& \limsup_{n \rightarrow \infty}  \frac{1}{n} \sum_{j=1}^{n} \mathbb{P}\Big(\frac{1}{a_n}T_j^{1 / \gamma_{n}} \leq y\Big)\\
&\leq& \limsup _{n \rightarrow \infty}  \frac{[nx]}{n}+ \limsup _{n
\rightarrow \infty}  \mathbb{P}\Big(\frac{1}{\gamma_{n}} \ln
\frac{T_{[nx]}}{\prod_{r=1}^{m_{n}}\Big(l_{r}+[nx]\Big)}-\ln F(x)
\leq  \ln \frac{y}{y+\delta}<0\Big)\\
&=&x=F^{*}(y+\delta).
\end{eqnarray*}
By setting $x=F^{*}(y-\delta)$ and $j=[nx]$ in \eqref{estimates} and
\eqref{eq3}, and using \eqref{limit3} again, we have
\begin{eqnarray*}
 \liminf _{n \rightarrow \infty}\frac{1}{n}\sum_{j=1}^{n} \mathbb{P}\Big(\frac{1}{a_n}T_j^{1 / \gamma_{n}} \leq y\big)
&\ge& \liminf _{n \rightarrow \infty}\frac{[nx]}{n}
\mathbb{P}\Big(\frac{1}{\gamma_{n}}
    \ln \frac{T_{[nx]}}{\prod_{r=1}^{m_{n}}\left(l_{r}+[nx]\right)}-\ln F(x)
    \leq  \ln \frac{y}{y-\delta}\Big)\\
    &=&x=F^{*}(y-\delta).
\end{eqnarray*}
Since $F^*(y)$ is continuous, we obtain \eqref{expectation1} by
letting $\delta$ tend to $0$.

Assume $y\leq F(0)$. For any $y_1 \in (F(0),1)$, we have
\[
\limsup_{n\to\infty}\frac{1}{n} \sum_{j=1}^{n}
\mathbb{P}\Big(\frac{1}{a_{n}} T_{j}^{1 / \gamma_{n}} \leqslant
y\Big)\leq \lim_{n\to\infty}\frac{1}{n}\sum_{j=1}^{n}
\mathbb{P}\Big(\frac{1}{a_{n}} T_{j}^{1 / \gamma_{n}} \leqslant
y_1\Big)= F^*(y_1),
\]
which tends to  $F^*(F(0))=0$ by letting $y_1\downarrow F(0)$ since
$F^*(y)$ is continuous. Similarly, when $y\geq 1$ we have for any
$y_2 \in (F(0),1)$
\[
\liminf_{n\to\infty}\frac{1}{n} \sum_{j=1}^{n}
\mathbb{P}\Big(\frac{1}{a_{n}} T_{j}^{1 / \gamma_{n}} \leqslant
y\Big)\ge \lim_{n\to\infty}\frac{1}{n}\sum_{j=1}^{n}
\mathbb{P}\Big(\frac{1}{a_{n}} T_{j}^{1 / \gamma_{n}} \leqslant
y_2\Big)= F^*(y_2),
\]
which tends to  $F^*(F(1))=1$ by letting $y_2\uparrow 1$.  In both
cases, \eqref{expectation1} still holds.

\vspace{10pt}

\noindent{\underline{\it Case 2}.} Assume $F$ is of Type II.

From Theorem~\ref{type}, we have
$\lim_{n\to\infty}\lambda_1(n)/\gamma_n=0$.   Since $\ln
G_n(x)=g_n(x)$ is bounded for any fixed $x\in (0,1)$, we have from
\eqref{limit1} that
\begin{equation}\label{wlln1}
\frac{1}{\gamma_n}\ln
\frac{T_{[\mathrm{nx}]}}{\prod_{r=1}^{m_n}(l_r+n)}\overset{p}\to 0
\end{equation}
for any $x\in (0,1]$.

Review $F^*$ in \eqref{case-b}.

When $y\le 0$,  \eqref{expectation1} is trivially true.

When $y\in (0,1)$, we have $\ln y<0$. For any $x\in (0,1)$, set
$j=[nx]$.  Then from \eqref{estimates}, \eqref{eq30} and
\eqref{wlln1} we have
\[
\limsup_{n \rightarrow \infty}  \frac{1}{n} \sum_{j=1}^{n}
\mathbb{P}\Big(\frac{1}{a_n}T_j^{1 / \gamma_{n}} \leq y\Big)\leq
\limsup _{n \rightarrow \infty}  \frac{[nx]}{n}+ \limsup _{n
\rightarrow \infty}  \mathbb{P}\Big(\frac{1}{\gamma_{n}} \ln
\frac{T_{[nx]}}{\prod_{r=1}^{m_{n}}\Big(l_{r}+[nx]\Big)}\le\ln
y\Big)=x.
\]
Since $x$ can be arbitrarily small, we have
\[
\lim_{n \rightarrow \infty}  \frac{1}{n} \sum_{j=1}^{n}
\mathbb{P}\Big(\frac{1}{a_n}T_j^{1 / \gamma_{n}} \leq
y\Big)=0=F^*(y).
\]
That is, \eqref{expectation1} is true.

When $y>1$, $\ln y>0$. Again,  for any $x\in (0,1)$, setting
$j=[nx]$ and using \eqref{estimates}, \eqref{eq30} and
\eqref{wlln1}, we have
\[
\liminf_{n\to\infty}\frac{1}{n} \sum_{j=1}^{n}
\mathbb{P}\Big(\frac{1}{a_{n}} T_{j}^{1 / \gamma_{n}} \leqslant
y\Big)\ge x.
\]
By letting $x\uparrow 1$ we have
\[
\lim_{n\to\infty}\frac{1}{n} \sum_{j=1}^{n}
\mathbb{P}\Big(\frac{1}{a_{n}} T_{j}^{1 / \gamma_{n}} \leqslant
y\Big)=1=F^*(y).
\]
This completes the proof of \eqref{expectation1}.

\noindent\underline{\textit{Case 3}}.  Assume $F$ is of Type III.

This time, we have $\lim_{n\to\infty}\lambda_1(n)/\gamma_n=\infty$.
\[
\frac{1}{\gamma_n}\ln
\frac{T_{[\mathrm{nx}]}}{\prod_{r=1}^{m_n}(l_r+n)}\overset{p}\to
-\infty
\]
for any $x\in (0,1)$. We can prove \eqref{expectation1} by using
similar lines to \textit{Case 2}. We omit the details.

The proof of Theorem~\ref{general} is completed.
\hfill$\blacksquare$

\vspace{10pt}

\noindent{\it Proof of Theorem~\ref{typeI}.} From
Theorems~\ref{general} and ~\ref{type}, we have  $\mu_n
\rightsquigarrow \mathrm{Unif}[0,2 \pi)\otimes \nu$ as $n\to\infty$,
where $\nu$ has a a density function
$f^*(r)=\frac{d}{dr}F^*(r)=\frac{1}{f(F^*(r))}$, $r\in [F(0),1]$.
Let $\Theta$ $R$ are two independent random variables, $\Theta$ is
uniformly distributed over $[0,2\pi)$ and $R$ has density function
$f^*$. Consider the transformation
$Z=R\exp(i\Theta)=R\cos(\Theta)+iR\sin(\Theta)=(R\cos(\Theta),
R\sin(\Theta))=:(X,Y)$. Note that the Jacobian for transformation
$(x,y)=(r\cos(\theta), r\sin(\theta))$ is $r=\sqrt{x^2+y^2}=|z|$,
where $z=re^{i\theta}=x+iy$. The joint density function of $Z=(X,Y)$
is given by $\frac{1}{2\pi}\frac{f^*(|z|)}{|z|}=\frac{1}{2\pi
f(F^*|z|)|z|}I(F(0)\le |z|\le 1)$.  Since $\hat\mu_n$ is obtained
under transformation $(x,y)=(r\cos(\theta), r\sin(\theta))$, by the
continuous mapping theorem, we $\hat\mu_n$ converges with
probability one to $Z=(X, Y)$ which has a joint density
$\frac{1}{2\pi f(F^*|z|)|z|}I(F(0)\le |z|\le 1)$.
\hfill$\blacksquare$

\vspace{10pt} \noindent{\it Proof of Theorem~\ref{typeII}}. Using
the same notations as in the proof for Theorem~\ref{typeI}, we have
$P(R=1)=1$. Therefore, we can easily conclude that $Z$ has a uniform
distribution on the unit circle.   \hfill$\blacksquare$

\vspace{20pt}

\noindent\textbf{Acknowledgements}. The authors would like to thank
an anonymous referee for his/her careful reading of the manuscript
and suggestion which has improved the layout of the manuscript.  The
research of Yongcheng Qi was supported in part by NSF Grant
DMS-1916014.

\baselineskip 12pt
\def\ref{\par\noindent\hangindent 25pt}

\end{document}